\documentclass[12pt]{article}
\usepackage{amsmath,amsthm,amssymb, amstext, amsopn, amsxtra}
\usepackage[all]{xy}
\usepackage{latexsym}
\usepackage{amsfonts}
\begin{document}
\newcommand{\nt}{\noindent}
\newcommand{\bs}{\bigskip}
\newcommand{\ms}{\medskip}
\newcommand{\mk}{\medskip}
\newcommand{\sk}{\smallskip}
\newcommand{\ep}{\varepsilon}
\newcommand{\m}{{\mathfrak m}}
\newcommand{\p}{{\mathfrak p}}
\newcommand{\dd}{\delta}
\newcommand{\A}{{\mathbb A}}
\newcommand{\R}{{\mathcal R}}
\newcommand{\M}{{\mathbb M}}
\newcommand{\ap}{\alpha^{\,\prime}}
\newcommand{\bp}{\beta^{\,\prime}}
\newcommand{\isom}{\cong}
\newcommand{\minus}{[\,e,e'\,]}
\newcommand{\plus}{\langle \,e,e'\,\rangle}



\[  \]

\begin{center}
\Large{\bf Commutators and Anti-Commutators of Idempotents in Rings}

\bigskip
\large{Dinesh Khurana and T.\,Y.~Lam} 

\end{center}

\begin{abstract}
\begin{small}
\nt We show that a ring $\,R\,$ has two idempotents $\,e,e'\,$
with an invertible commutator $\,ee'-e'e\,$ if and only if $\,R
\cong {\mathbb M}_2(S)\,$ for a ring $\,S\,$ in which $\,1\,$
is a sum of two units.  In this case, the ``anti-commutator"
$\,ee'+e'e\,$ is automatically invertible, so we study
also the broader class of rings having such an invertible
anti-commutator.  Simple artinian rings $\,R\,$ (along with 
other related classes of matrix rings) with one of the above 
properties are completely determined. In this study, we
also arrive at various new criteria for {\it general\/}
$\,2\times 2\,$ matrix rings. For instance, $R\,$ is such a
matrix ring if and only if it has an invertible commutator
$\,er-re\,$ where $\,e^2=e$.
\end{small}
\end{abstract}

\mk\nt 
{\bf \S1. \ Introduction}

\bs
The work in this paper was inspired by an insightful exercise 
of Kaplansky in his 1968 book [Ka] on rings of operators. 
Given two idempotents $\,e,e'\,$ in a ring $\,R\,$ with the 
property that their commutator $\,ee'-e'e\,$ is invertible, 
``Exercise 6'' in [Ka:~p.\,25] asked the readers to show that 
the idempotents $\,e,e',1-e,1-e'\,$ are pairwise isomorphic; 
in other words, the four principal right ideals generated
by them are isomorphic as right\break
$R$-modules. This exercise was intended for prospective
use toward the proofs of Theorem 60 and Theorem 61 on Baer
$\,\ast\,$-rings in [Ka:~pp.\,91--96]. For his Exercise 6
on p.\,25, Kaplansky kindly offered his readers a ``Hint''.
However, rendered with the author's trademark brevity,
this ``Hint'' was itself no less than another substantial
exercise. In all fairness, a full solution of Kaplansky's
``Exercise 6'' proved to be a considerable challenge by
any yardstick.  

\bs
That was the situation about fifty years ago, in 1967--68. 
Nowadays, with so much more known about idempotents and idempotent
identities in rings, a rather natural new solution can be given
for Kaplansky's ``Exercise 6''. More remarkably, the theme of this
exercise can be further developed so as to give a full-fledged
characterization theorem for the rings $\,R\,$ that appeared in 
the exercise. To explain this from a more general point of view,
we'll use the notation $\,[x,y]\,$ for the {\it commutator\/}
$\,xy-yx$, and the notation $\,\langle x,y\rangle\,$ for the
{\it anti-commutator\/} $\,xy+yx$, for any two elements $\,x,y\,$
in a ring $\,R$. For the main purposes of this paper, it is
convenient to introduce the following two ring-theoretic
properties, where $\,{\rm idem}\,(R)\,$ will henceforth denote
the set of idempotents in the ring $\,R$, and $\,{\rm U}(R)\,$
will denote the group of units of $\,R$.

\bs\nt
{\bf Property K\,:} {\it There exist $\,e,e'\in {\rm idem}\,(R)\,$
such that $\,\minus \in {\rm U}(R)$.}

\mk\nt
{\bf Property $\boldsymbol{\overline{{\rm K}}}\,$:} {\it There exist
$\,e,e'\in {\rm idem}\,(R)\,$ such that $\,\plus \in {\rm U}(R)$.}

\bs
The easiest class of rings to deal with in the investigation
of these properties is the class of {\it abelian rings\/};
that is, rings in which all idempotents are central.  It is
easy to see that such a ring $\,R\,$ has Property K iff
$\,R=0$, while it has Property $\overline{{\rm K}}\,$
iff $\,2\in {\rm U}(R)$. (In fact, in view of the fact that
$\,\langle\,1,1\,\rangle =2$, any ring with $\,2\in
{\rm U}(R)\,$ has Property $\overline{{\rm K}}$.) Another
easy observation we can make is that, if $\,R\to R'\,$ is
a (unital) ring homomorphism, then $\,R\,$ having one of
the two properties above implies that $\,R'\,$ has the same
property, whereby any matrix ring $\,{\mathbb M}_n(R)\,$
will also have the same property.

\bs
The first two main results in this paper are Theorem A and
Theorem B below, where the former (to be proved in \S2)
explains our choice of the notation $\overline{{\rm K}}\,$
since it has the obvious consequence that the class of rings
with Property K is (properly) contained in the class of
rings with Property $\overline{{\rm K}}$. 

\bs\nt
{\bf Theorem A.} {\it For any $\,e,e'\in {\rm idem}\,(R)$,
$\,\minus \in {\rm U}(R)\,$ implies that} $\,\plus \in
{\rm U}(R)\,$ ({\it though not conversely\/}). {\it In
particular, if $\,R\,$ has {\rm Property K\/}, then it
has {\rm Property\/} $\overline{{\rm K}}$.} 

\bs\nt
{\bf Theorem B.} {\it A ring $\,R\,$ has {\rm Property K\/} 
iff $\,R\cong {\mathbb M}_2(S)\,$ for a ring $\,S\,$
in which $\,1\in {\rm U}(S)+{\rm U}(S)$.}

\bs
From our more general perspective, the solution of Kaplansky's
Exercise 6 is just a part of the work needed for proving the
``only if'' part of Theorem B. In \S3, this theorem is proved
by using judiciously various idempotent identities of Nicholson
[Ni], Kato [Kt], and Koliha-Rako\v{c}evi\'{c} [KR$_1$, KR$_2$],
which we first develop {\it ad hoc\/} in \S2.  The classical
notion of Bott-Duffin invertibility (relative to an idempotent)
introduced in [BD] turns out to play a significant role in
proving Theorem B, so the basic ingredients of the Bott-Duffin
theory are briefly recalled at the beginning of \S2 as well.

\bs
With Theorem A and Theorem B at our disposal, a natural question
to ask is when would a matrix ring $\,{\mathbb M}_n(T)\,$ (over
a given type of base rings $\,T\,$) have Property K or Property
$\overline{{\rm K}}$.  For instance, in the case where $\,T\,$
is a division ring, this would amount to asking which simple
artinian rings would have Property K or Property
$\overline{{\rm K}}$. Indeed, taking $\,T\,$ to be a local
ring or a nonzero commutative ring, we have the following
result in \S3 on matrix rings $\,{\mathbb M}_n(T)\,$ with
Property K.

\bs\nt
{\bf Theorem C.}  {\it Let $\,R={\mathbb M}_n(T)$, where 
$\,T\,$ is a local ring or a nonzero commutative ring.
Then $\,R\,$ has {\rm Property K\/} iff either
$\,n\in \{4,6,8,\dots\,\}$, or $\,n=2\,$ and $\,1_T\in
{\rm U}(T)+{\rm U}(T)$.} ({\it Note that in the case where
$\,(T, {\mathfrak m})\,$ is local ring, the condition $\,1_T
\in {\rm U}(T)+{\rm U}(T)\,$ would amount to the simpler
statement that $\,|\,T/{\mathfrak m}\,|>2$.})

\bs
As for Property $\overline{{\rm K}}$, the corresponding
result (to be proved in \S4) is the following.

\bs\nt
{\bf Theorem D.} {\it Let $\,T\,$ be a local ring or
a commutative ring. Then $\,R={\mathbb M}_n(T)\,$ has
{\rm Property $\overline{{\rm K}}$\/} iff we are in
one of the following cases\/}:

\mk\nt
(1) $\,n\in \{4,6,8,\dots \,\}$.\\
(2) $\,n=2\,$ {\it and $\,1_T\in {\rm U}(T)+{\rm U}(T)$.} \\
(3) $\,n\,$ {\it is odd and $\,2\in {\rm U}(T)$.}

\bs
As a spin-off of our investigations on invertible commutators
of the form $\,\minus\,$ where $\,e,e'\in {\rm idem}\,(R)$,
we also consider the case where $\,e'\,$ is allowed to be
an arbitrary element of $\,R$, while $\,e\,$ is replaced by
an element with some specific property ``comparable'' to
idempotency.  Working with these more general commutators (and
anti-commutators) and building on the work of Fuchs-Maxson-Pilz
in [FMP:~(III.2)], we obtain in \S5 a number of new criteria
for $\,2\times 2\,$ matrix rings, three of which are summarized
as follows.

\bs\nt
{\bf Theorem E.} {\it A ring $\,R\,$ is a $\,2\times 2\,$ 
matrix ring over some other ring iff there is an invertible
commutator $\,[\,e,\,\ast\,]\,$ where $\,e^2=e$, iff there
is an invertible commutator $\,[\,p,\,\ast\,]\,$ where
$\,p^2=0$, iff there is an invertible commutator
$\,[\,u,\,\ast\,]\,$ where $\,u^2=1$.} 

\bs
The terminology and notations introduced so far in this 
Introduction will be used freely throughout the paper. For any
ring $\,R$, $\,{\rm rad}\,(R)\,$ denotes the Jacobson radical 
of $\,R$, and the words ``exchange ring'' will be used in
the sense of Warfield [Wa] and Nicholson [Ni]. By saying that
a matrix $\,M\in {\mathbb M}_n(T)\,$ is {\it diagonalizable,}
we'll mean that $\,M\,$ is similar to a diagonal matrix in
$\,{\mathbb M}_n(T)$.  Other standard terminology and
conventions in ring theory follow mainly those in [Go],
[La$_1$], and [La$_3$]. 

\mk

\bs\nt 
{\bf \S2. \ Idempotent Identities of Kato and Koliha-Rako\v{c}evi\'c}

\bs
One key ingredient used in this beginning section is the basic
notion of Bott-Duffin invertibility introduced in the early
paper [BD]. In order to give a relatively self-contained exposition
of our results, we'll start by recalling some definitions and
facts from [BD]. For any idempotent $\,e\,$ in a ring $\,R$, an
element $\,a\in R\,$ is said to be {\it Bott-Duffin invertible\/}
relative to $\,e\,$ if $\,eae\in {\rm U}(eRe)$. In this case, the
inverse of $\,eae\,$ in the corner ring $\,eRe\,$ is said to be
the Bott-Duffin inverse of $\,a\,$ relative to $\,e$. For the
sake of completeness, we state and prove the following classical
characterization result for Bott-Duffin invertibility.

\bs\nt
{\bf Theorem 2.1.} {\it For $\,a\in R\,$ and $\,e\in {\rm idem}\,(R)$, 
the following are equivalent\/}:

\mk\nt
(1) {\it $\,a\,$ is Bott-Duffin invertible relative to $\,e$.}\\
(2) $\,1-e+ae\in {\rm U}(R)$. \;\;
(3) $\,1-e+ea\in {\rm U}(R)$. \;\;
(4) $\,1-e+eae\in {\rm U}(R)$.

\bs\nt
{\bf Proof.} For any ring $\,R$, Jacobson's Lemma (see, 
e.g.~[La$_2$:~Exercise 1.6]) states that, for any $\,x,y
\in R$, $\,1-xy\in {\rm U}(R)\,$ iff $\,1-yx\in {\rm U}(R)$.  
Since $\,1-e+ae=1-(1-ae)\,e$, this Lemma gives
$\,(2)\Leftrightarrow (4)$.  By left-right symmetry,
we have also $\,(3)\Leftrightarrow (4)$.  To see that
$\,(4)\Leftrightarrow (1)$, let $\,f=1-e$.  With respect to
the idempotent $\,e$, the element $\,f+eae\,$ has a diagonal
Peirce decomposition matrix $\footnotesize{\begin{pmatrix}
eae&0\\0&f\end{pmatrix}}$. Such a matrix is invertible iff
$\,eae\in {\rm U}(eR\,e)$, which is the defining condition
for (1).  This shows that $\,(4)\Leftrightarrow (1)$. \qed

\bs
Next, we state the following key result from Kato's book 
[Kt:~(I.4.34), (I.4.44)].

\bs\nt
{\bf Proposition 2.2. (Kato's Identities)}
{\it For any $\,e,e'\in {\rm idem}\,(R)$, we have
$$
(e-e')^2+(1-e-e')^2=1, \leqno (2.3)
$$
which amounts to 
$$
(e-e')^2=(e+e')\,(2-e-e'), \leqno (2.4)
$$ 
where the two factors on the {\rm RHS\/} commute. Also,
writing $\,r=1-e+e'e\,$ and $\,s=1-e'+ee'$, we have}
$$ 
rs = sr = (1+e-e')\,(1-e+e') = (1-e-e')^2. \leqno (2.5)
$$
\sk\nt
{\bf Proof.} The identity (2.3) is easily verified by a direct
expansion of the LHS. By transposing the term $\,(1-e-e')^2\,$
to the RHS and using the factorization $\,1-x^2=(1-x)\,(1+x)$, 
we get the identity (2.4). Similarly, by transposing the term
$\,(e-e')^2\,$ in (2.3) to the RHS and using the factorization
$\,1-y^2=(1+y)\,(1-y)$, we get the last equality in (2.5).
The first two equalities in (2.5) are easily verified by
direct expansions.\qed 

\bs\nt
{\bf Corollary 2.6.} {\it Two idempotents $\,e,e'\in R\,$ are
Bott-Duffin invertible relative to each other iff 
$\,1-e-e'\in {\rm U}(R)$.}

\bs\nt
{\bf Proof.} This follows from (2.5), applied in conjunction 
with Theorem 2.1.\qed

\bs\nt
{\bf Remark 2.7.} In general, $\,e'\,$ being Bott-Duffin invertible
relative to $\,e\,$ alone does not imply that $\,e\,$ is
Bott-Duffin invertible relative to $\,e'$. For instance, if $
\,e'=1\neq e$, then $\,e'\,$ is Bott-Duffin invertible relative
to $\,e$, but $\,e\,$ is not Bott-Duffin invertible relative to
$\,e'$.  In this example, we have $\,1-e-e'=-e\notin {\rm U}(R)$.
    
\bs
In [KR$_1$:~Theorem 3.5], it was proved that, for any idempotents
$\,e,e'\,$ in any ring, $\,e-e'\in {\rm U}(R)\,$ iff $\,e+e'\in 
{\rm U}(R)\,$ and $\,1-ee'\in {\rm U}(R)$. (The hypothesis $\,2\in 
{\rm U}(R)\,$ was included in the statement of this theorem, but
the result can be seen to be true without such a hypothesis.)
Using here the crucial identity (2.4), we'll prove the following 
closely related result.

\bs\nt
{\bf Theorem 2.8.} {\it For any $\,e,e'\in {\rm idem}\,(R)$,
we have $\,e-e'\in {\rm U}(R)\,$ iff $\,e+e'\in {\rm U}(R)\,$ 
and $\,f+f'\in {\rm U}(R)$, where $\,f=1-e\,$ and $\,f'=1-e'$.
In this case, $\,e\,$ is similar to $\,1-e'$.}

\bs\nt
{\bf Proof.} The ``iff" statement follows from the identity (2.4)
(since $\,2-e-e'=f+f'$). The last statement is a special case 
of a result of Nicholson [Ni:~Proposition 1.8] on the clean 
representations of elements in a ring. To make our exposition
self-contained, we recall Nicholson's proof (in our special
case). If $\,u:=e-e'\in {\rm U}(R)$, then
$$
(1-e')\,u=(1-e')\,(e-e')=(1-e')\,e=(e-e')\,e=ue. 
$$
Thus, $\,e=u^{-1}(1-e')\,u\,$ is similar to $\,1-e'$.\qed

\bs\nt
{\bf Remark 2.9.}  In general, $\,e+e'\in {\rm U}(R)\,$ alone
need not imply that $\,e-e'\in {\rm U}(R)$.  For instance, in
any nonzero ring $\,R\,$ in which $\,2\in {\rm U}(R)$, the
idempotents $\,e=e'=1\,$ have the property that $\,e+e'=2\in
{\rm U}(R)$, but $\,e-e'=0\notin {\rm U}(R)$.

\bs
Next we give a systematic derivation for two useful identities 
on the commutator $\,\minus =ee'-e'e\,$ and the anti-commutator
$\,\plus =ee'+e'e\,$ from [KR$_2$:~p.\,103] and [KR$_1$:~p.\,289].  

\bs\nt
{\bf Proposition 2.10. (Koliha-Rako\v{c}evi\'c Identities)}
{\it For any $\,e,e'\in {\rm idem}\,(R)\,$ and $f=1-e$, 
we have\/}
$$
\minus = (e-e')\,(e'-f) = (f-e')\,(e-e'); \leqno (2.11)
$$
\vspace{-0.4in}
$$
\plus = (e+e')\,(e'-f) = (e'-f)\,(e+e').\leqno (2.12)
$$

\sk\nt
{\bf Proof.} From $\,(e+e')^2=e+e'+\plus $, transposition
of the term $\,e+e'\,$ gives the two equalities in (2.12). 
Similarly, from $\,(e-e')\,(e+e')=e-e'+\minus$, transposition
of the term $\,e-e'\,$ gives the first equality in (2.11),
and the second equality follows by working instead with
$\,(e+e')\,(e-e')=e-e'-\minus $.\qed

\bs
Kaplansky asserted in [Ka:~p.\,25, Exercise 6] that, under
the assumption $\,\minus \in {\rm U}(R)$, the idempotents
$\,e,e'\in R\,$ are Bott-Duffin invertible relative to each
other.  With the help of Proposition 2.10, one can prove the
more precise result in (1) below, and by the same token,
prove also its complete analogue in (2) for anti-commutators.   

\bs\nt
{\bf Theorem 2.13.} {\it For any $\,e,e'\in {\rm idem}\,(R)\,$
and $\,f=1-e$, the following hold.}

\mk\nt
(1) {\it $\,\minus \in {\rm U}(R)\,$ iff $\,e-e'\in {\rm U}(R)\,$ 
and $\,f-e'\in {\rm U}(R)$, iff $\,e-e'\in {\rm U}(R)\,$ and
$\,e,e'\,$ are Bott-Duffin invertible relative to each other.}

\sk\nt
(2) {\it $\,\plus \in {\rm U}(R)\,$ iff $\,e+e'\in {\rm U}(R)\,$ 
and $\,f-e'\in {\rm U}(R)$, iff $\,e+e'\in {\rm U}(R)\,$ and 
$\,e,e'\,$ are Bott-Duffin invertible relative to each other.}

\sk\nt
(3) {\it $\,\minus \in {\rm U}(R)\,$ implies that 
$\,\plus \in {\rm U}(R)$.}

\bs\nt
{\bf Proof.} For (1), the first ``iff" statement follows from 
(2.11); see [KR$_2$:~(3.6)].  The second ``iff" statement then 
follows from Corollary 2.6.  (2) is proved similarly by using 
(2.11) and Corollary 2.6 (although the first ``iff" statement
in (2) was stated in [KR$_1$:~Theorem 3.6] with an extra 
assumption that $\,2\in {\rm U}(R)$).  Finally, (3) follows 
from (1) and (2), since $\,e-e'\in {\rm U}(R)\,\Rightarrow\,
e+e'\in {\rm U}(R)\,$ according to Theorem 2.8. \qed

\bs
Recalling from \S1 that ``Property K'' (respectively, ``Property
$\overline{{\rm K}}\,$'') on a ring $\,R\,$ means the existence
of $\,e,e' \in {\rm idem}\,(R)\,$ such that $\,\minus \in
{\rm U}(R)\,$ (respectively, $\,\plus \in {\rm U}(R)$), we can
draw the following conclusions from Theorem 2.13.

\bs\nt
{\bf Corollary 2.14.} {\it If a ring $\,R\,$ has
{\rm Property K}, then it has {\rm Property $\overline{{\rm K}}$.}
However, the converse of this statement fails in general.}

\bs\nt
{\bf Proof.} The first statement here follows from part
(3) of Theorem 2.13. The second statement follows from the
observation (made in \S1) that any ring $\,R\,$ with $\,2
\in {\rm U}(R)\,$ has Property $\overline{{\rm K}}$, while
an abelian ring $\,R\neq 0\,$ cannot have Property K. \qed

\bs
For future reference, we record in the following a useful
result on the Jacobson radical and an immediate consequence
thereof.

\bs\nt
{\bf Proposition 2.15.} {\it Let $\,J\subseteq {\rm rad}\,(R)\,$
be an ideal of $\,R\,$ such that idempotents lift modulo $\,J$.
Then $\,R\,$ has {\rm Property K} iff $\,R/J\,$ does. The same
statement holds for {\rm Property $\overline{{\rm K}}$.}}

\bs\nt
{\bf Proof.} As we have observed in \S1, the ``only if'' part
holds without any assumptions on the ideal $\,J$. In the case
where idempotents lift modulo $\,J\subseteq {\rm rad}\,(R)$,
the ``if'' part follows from the observation that an element
$\,s\in R\,$ is a unit iff $\,\overline{s}\,$ is a unit in
$\,\overline{R}=R/J$.\qed.

\bs\nt
{\bf Corollary 2.16.} {\it Let $\,J\,$ be an ideal of $\,R\,$
such that one of the following holds\/}:

\mk\nt
(1) $\,J\,$ {\it is a nil ideal.}\\
(2) $\,J\subseteq {\rm rad}\,(R)\,$ {\it and $\,R\,$ is an
exchange ring.}

\mk\nt
{\it Then the conclusions of {\rm Proposition 2.15\/} hold
for the two rings $\,R\,$ and $\,R/J$.} 

\bs\nt
{\bf Proof.} (1) If $\,J\,$ is a nil ideal, we have $\,J
\subseteq {\rm rad}\,(R)\,$ by [La$_1$:~(4.11)]. On the other
hand, by [La$_1$:~(21.28)], idempotents lift modulo $\,J$.
Thus, Proposition 2.15 applies.

\mk\nt
(2) If $\,R\,$ is an exchange ring, idempotents can be lifted
modulo any ideal in $\,R\,$ (according to [Ni]), so again
Proposition 2.15 applies.\qed

\bs
We shall conclude the discussions in this section by using
Corollary 2.16 to prove the following two easy results on
upper triangular rings.

\bs\nt
{\bf Proposition 2.17.} {\it Let $\,S,\,T\,$ be two rings, and
let $\,M\,$ be any $\,(S,T)$-bimodule. Then the formal triangular
ring $\,R=\footnotesize{\begin{pmatrix}S&M\\0&T\end{pmatrix}}$
has} {\rm Property K} ({\it respectively,} Property
$\overline{{\rm K}}\,$) {\it iff $\,S\,$ and $\,T\,$ both do.}

\bs\nt
{\bf Proof.} The ``only if'' part is clear, since $\,S\,$ and
$\,T\,$ can both be thought of as factor rings of $\,R$.
Conversely, assume that $\,S\,$ and $\,T\,$ have Property K.
Then $\,S\times T\,$ also does. As
$\,J\!:=\footnotesize{\begin{pmatrix}0&M\\0&0\end{pmatrix}}$
is a square-zero ideal in $\,R$, and $\,R/J\cong S\times T$,
it follows from Case (1) of Corollary 2.16 that $\,R\,$ also
has Property K. The case where $\,S\,$ and $\,T\,$ both have
Property $\overline{{\rm K}}\,$ is similar.\qed

\bs
In the case of $\,n\times n\,$ upper triangular matrices
over a given ring $\,S$, the ideal $\,J\,$ of strictly upper
triangular matrices has the property that $\,J^{\,n}=0$.
Thus, the same method of proof (using Case (1) of Corollary
2.16) yields the following similar result.

\bs\nt
{\bf Proposition 2.18.} {\it Let $\,R={\mathbb T}_n(S)\,$ be the
ring of $\,n\times n\,$ upper triangular matrices over a ring
$\,S$. Then $\,R\,$ has} {\rm Property K} ({\it respectively,}
Property $\overline{{\rm K}}\,$) {\it iff $\,S\,$ does.}

\mk

\bs\nt
{\bf \S3. \ Rings with Property K}

\bs 
To study rings with the Property K, we begin by giving a 
streamlined proof for the rest of Kaplansky's ``Exercise 6" 
from [Ka:~p.\,25]; see also Theorem 60 in [Ka:~p.\,96].  The 
proof to be given below is in some sense much easier to grasp
than that given by Kaplansky since it makes full use of the
efficient tools developed in \S2.

\bs\nt
{\bf Theorem 3.1.} {\it If $\,e,e'\in {\rm idem}\,(R)\,$ are 
such that $\,\minus \in {\rm U}(R)$, then the following hold.}

\mk\nt
(1) {\it The four idempotents $\,e,e',1-e,1-e'\,$ are pairwise 
similar.}\\
(2) $\,R=eR\oplus e'R=R\,e\oplus R\,e'$.

\bs\nt
{\bf Proof.} (1) Let $\,f=1-e$. By Theorem 2.13(1), we have $\,e-e'
\in {\rm U}(R)\,$ and $\,f-e'\in {\rm U}(R)$. Therefore, by the
last part of Theorem 2.8, $\,e'\,$ is similar to $\,f\,$ as well 
as to $\,e$. By symmetry, it follows that $\,e,e',1-e,1-e'\,$ 
are pairwise similar.

\mk\nt
(2) By symmetry, it suffices to prove that $\,R=eR\oplus e'R$.
Since $\,e-e'\in {\rm U}(R)$, clearly $\,eR+e'R=R$.  Also, 
for any $\,x\in eR\cap e'R$, we have $\,x=ex=e'x$.  Thus, 
$\,(e-e')\,x=0$, and hence $\,x=0$. This proves that
$\,R=eR\oplus e'R$. (In conclusion, we note that (2) can also 
be deduced as a consequence of [KR$_2$:~Theorem 3.2].)\qed

\bs\nt
{\bf Remark 3.2.} In general, if two idempotents $\,e,e'\,$
in a ring $\,R\,$ are such that $\,e,e',\break
1-e,1-e'\,$ are pairwise similar, the commutator $\,\minus \,$
need not be a unit in $R$, and $\,e'\,$ need not be
Bott-Duffin invertible relative to $\,e$.  For instance, if
$\,R={\mathbb M}_2(S)\,$ over a nonzero ring $\,S\,$ and
$\,E_{ij}\,$ are the matrix units, then $\,e=E_{11}\,$ and
$\,e'=E_{22}\,$ certainly have the similarity properties
mentioned above, but $\,\minus =0$, and $\,ee'e=0\,$ (so
$\,e'\,$ is not Bott-Duffin invertible relative to $\,e$).
On the other hand, if we have only $\,[\,e,r]\in {\rm U}(R)\,$
where $\,r\in R\,$ (or even $\,r\in {\rm U}(R)$), it also
does not follow that $\,r\,$ is Bott-Duffin invertible
relative to $\,e$, as is shown by the example $\,e=E_{11}\,$
and $\,r=E_{12}+E_{21}$, for which $\,[\,e,r]=E_{12}-E_{21}
\in {\rm U}(R)$, but again $\,ere=0$.

\bs 
Now we are in a position to prove Theorem B stated in the
Introduction. The remarkable thing about this theorem is
that a {\it single\/} commutator property on a ring turns
out to be enough to determine the structure of the whole
ring.  Just for the record, we note that it was Kaplansky's
``Exercise 6'' in [Ka:~p.\,25] which had given us the impetus
to prove the ``only if'' part of this theorem.

\bs\nt
{\bf Theorem 3.3.}  {\it A ring $\,R\,$ has {\rm Property K\/}
iff $\,R\cong {\mathbb M}_2(S)\,$ for a ring $\,S\,$ in which 
$\,1\in {\rm U}(S)+{\rm U}(S)$.}

\bs\nt
{\bf Proof.} First assume that $\,1=a+b\in S\,$ where
$\,a,b\in {\rm U}(S)$. In the matrix ring $\,{\mathbb M}_2(S)$,
we check easily that $\,e=\footnotesize{\begin{pmatrix}a&b\\a&b
\end{pmatrix}}$ is an idempotent, and that the commutator
$\,[\,e,E_{11}\,]=\footnotesize{\begin{pmatrix}0&-b\\a&0
\end{pmatrix}}$ is a unit. Thus, $\,{\mathbb M}_2(S)\,$ has
Property K. (In testimony to the implication statement in Theorem
2.13(3), the anti-commutator $\,\langle e, E_{11}\,\rangle
=\footnotesize{\begin{pmatrix}2\,a&b\\a&0\end{pmatrix}}$
is indeed a unit too, so the ring $\,{\mathbb M}_2(S)\,$ also
has Property $\overline{{\rm K}}$ as expected.)

\mk
Conversely, if a ring $\,R\,$ has Property K, fix two idempotents
$\,e,e'\in R\,$ such that $\,\minus \in {\rm U}(R)$. By Theorem
3.1, $\,f\!:=1-e\,$ is similar to $\,e$. This implies that
$\,fR\cong eR$, and so $R_R=eR\oplus fR \cong eR\oplus eR$.
Taking endomorphism rings, we have
$$
R \cong {\rm End}_R(eR\oplus eR) \cong {\mathbb M}_2(S) 
$$
where $\,S\!:={\rm End}_R(eR)\,$ can be identified with the
corner ring $\,eR\,e$. By Theorem 2.13(1), $\,ee'e\in {\rm U}(S)$.
Similarly, for $\,f'\!:=1-e'$, the fact that $\,[\,e,f'\,]=-\minus
\in {\rm U}(R)\,$ implies that $\,ef'e\in {\rm U}(S)\,$ too. Thus,
$\,1_S=e=ee'e+ef'e\in {\rm U}(S)+{\rm U}(S)$, as desired.\qed

\bs
In connection with the condition $\,1\in {\rm U}(S)+{\rm U}(S)\,$
appearing in Theorem 3.3, we should point out that it is in fact
equivalent to the following ostensibly stronger condition: for
any integer $\,n\geq 1$, any unit in $\,S\,$ is a sum of $\,n\,$
units in $\,S$. In the case where $\,S\,$ is a unit-regular ring,
the condition $\,1_S\in {\rm U}(S)+{\rm U}(S)\,$ can also be
self-strengthened into $\,S={\rm U}(S)+{\rm U}(S)$, according
to [GW:~Theorem 3.8].  In the following, we will state two
interesting consequences of Theorem 3.3.

\bs\nt
{\bf Corollary 3.4.} {\it A ring $\,R\,$ has {\rm Property K\/}
with $\,2\in {\rm U}(R)\,$ iff $\,R\cong {\mathbb M}_2(S)\,$
for some ring $\,S\,$ with $\,2\in {\rm U}(S)$.}

\bs\nt
{\bf Proof.} The ``if'' part follows from Theorem 3.3 since
$\,2\in {\rm U}(S)\,$ implies that $\,1=2^{-1}+2^{-1}\in
{\rm U}(S)+{\rm U}(S)$, as well as $\,2\in
{\rm U}\bigl({\mathbb M}_2(S)\bigr)$. Conversely, if a ring
$\,R\,$ has {\rm Property K\/} along with $\,2\in {\rm U}(R)$,
Theorem 3.3 implies that $\,R\cong {\mathbb M}_2(S)\,$ for
some ring $\,S$. Clearly, the fact that $\,2\in {\rm U}(R)\,$
implies that $\,2\in {\rm U}(S)$.\qed

\bs
A second application of the ``if'' part of Theorem 3.3 is
that it gives a good supply of examples of matrix rings
having the Property K, as follows.

\bs\nt
{\bf Corollary 3.5.} {\it If $\,R={\mathbb M}_{2m}(T)\,$}
({\it for any ring $\,T$\/}) {\it where $\,m\geq 2$, or 
$\,R={\rm End}_D(V)\,$ for some infinite-dimensional right
vector space $\,V$ over a division ring $\,D$, then $\,R\,$ 
has\/} {\rm Property K.}

\bs\nt
{\bf Proof.} First assume that $\,R={\mathbb M}_{2m}(T)\,$
where $\,m\geq 2$. By a result of Henriksen [He], the identity
matrix $\,I_m\,$ is a sum of two units in $\,{\mathbb M}_m(T)$.
Thus, by Theorem 3.3, $\,R\cong {\mathbb M}_{2}
\bigl({\mathbb M}_m(T)\bigr)\,$ has Property K. Finally,
assume that $\,R={\rm End}\,(V_D)\,$ as in the statement of
the Corollary.  By [La$_3$:~Example 1.4], we have $\,R\cong
{\mathbb M}_n(R)\,$ for every $\,n\geq 1$. Applying this for
$\,n=4\,$ (for instance), we are back to the case treated
above. \qed

\bs\nt
{\bf Example 3.6.} In general, if $\,T\,$ is a ring such that
$\,R\!:={\mathbb M}_2(T)\,$ has Property K, Theorem 3.3 only
says that $\,R\cong {\mathbb M}_2(S)\,$ for some ring $\,S\,$
in which $\,1\,$ is a sum of two units.  Since $\,T\,$ may
not be isomorphic to $\,S$, this {\it does not\/} imply that
$\,1_T\in {\rm U}(T)+{\rm U}(T)$.  An explicit example to
illustrate the possible failure of $\,1_T\in {\rm U}(T)
+{\rm U}(T)\,$ is as follows.  For any field $\,k$, we
construct after Leavitt [Le] and Cohn [Co$_1$] a $\,k$-algebra
$\,T\,$ with a generic invertible $\,3\times 2\,$ matrix.
This matrix defines a right $\,T$-module isomorphism $\,T^2
\cong T^3$, which induces another isomorphism $\,T^2\cong T^4$.
Taking $\,T$-endomorphism rings gives a ring isomorphism
$\,R\!:={\mathbb M}_2(T)\cong {\mathbb M}_4(T)$. Thus, it
follows from Corollary 3.5 that $\,R\,$ has Property K. By
a recent result of G.~Bergman [Be$_2$] based on his earlier
work [Be$_1$], $\,{\rm U}(T)={\rm U}(k)$. In particular, if
we take $\,k\,$ to be the field of two elements, we'll have
$\,{\rm U}(T)=\{1\}$, in which case clearly $\,1_T\notin
{\rm U}(T)+{\rm U}(T)$.

\bs
Example 3.6 leads naturally to the following

\bs\nt
{\bf Question 3.7.} {\it If $\,T\,$ is a ring such that
$\,{\mathbb M}_2(T)\,$ has {\rm Property K}, under what
additional assumptions on $\,T\,$ can we conclude that} 
$\,1_T \in {\rm U}(T)+{\rm U}(T)\,$?

\bs
The difficulty in dealing with this question stems mainly from
the fact that $\,1_T \in {\rm U}(T)+{\rm U}(T)\,$ is {\it not\/}
a Morita invariant property of rings. Nevertheless, there are many
classes of rings $\,T\,$ for which we can answer Question 3.7
in a satisfactory way.  First, recall from [La$_3$:~\S17C] that
a ring $\,T\,$ is said to be
$\,\boldsymbol{{\mathbb M}_n}$-{\bf unique\/} if, for any ring
$\,S$, $\,{\mathbb M}_n(T)\cong {\mathbb M}_n(S)\,$ implies
that $\,T\cong S$. Taking stock in this definition, we do have
the following partial positive answer to Question 3.7.

\bs\nt
{\bf Theorem 3.8.} {\it Let $\,R={\mathbb M}_2(T)\,$ where
$\,T\,$ satisfies one of the following conditions\/}:

\sk\nt
(1) $\,T\,$ {\it is $\,{\mathbb M}_2$-unique.}\\
(2) $\,T\,$ {\it is an abelian ring and every idempotent in
$\,R\,$ is diagonalizable.}

\sk\nt
{\it Then $\,R\,$ has {\rm Property K\/} iff $\,1_T\in
{\rm U}(T)+{\rm U}(T)$.}

\bs\nt
{\bf Proof.} The ``if'' part is true without any assumptions
on $\,T\,$ by Theorem 3.3. The ``only if'' part in the two
cases (1) and (2) will be handled separately, as follows.

\mk\nt
{\bf Case (1).} If $\,R\,$ has Property K, Theorem 3.3
implies that $\,R \cong {\mathbb M}_2(S)\,$ for some ring
$\,S\,$ with $\,1_S\in {\rm U}(S)+{\rm U}(S)$. Given the
$\,{\mathbb M}_2$-unique assumption in Case (1), we have
$\,T\cong S$, and so $\,1_T\in {\rm U}(T)+{\rm U}(T)$.

\mk\nt
{\bf Case (2).} If $\,R\,$ has Property K, fix a commutator
$\,\minus \in {\rm U}(R)\,$ with $\,e,e'\in {\rm idem}\,(R)$.
Given the hypothesis in this case, we may assume (after a
conjugation) that $\,e'={\rm diag}\,(s,t)$, where $\,s,t
\in T\,$ are necessarily central idempotents. Writing 
$\,e=\footnotesize{\begin{pmatrix}a&b\\c&d\end{pmatrix}}$,
we have $\,\minus =\footnotesize{\begin{pmatrix}0&b\,(t-s)\\ 
c\,(s-t)&0\end{pmatrix}}\in {\rm U}(R)$, which implies that 
$\,b,c\in {\rm U}(T)$. From the equation $\,e^2=e$, we get 
$\,a^2+b\,c=a$, so $\,a,\,1-a\,$ are units in $\,T$,
with sum $1$.\qed

\bs
Having proved Theorem 3.8, we will mention in
(3.9)\,--\,(3.13) below some of the more important classes
of base rings $\,T\,$ to which the theorem can be applied.

\bs\nt
{\bf (3.9)} The most obvious class of examples of a ring
$\,T\,$ satisfying the hypothesis (2) of Theorem 3.8 is
given by the projective-free rings defined by P.\,M.~Cohn in
[Co$_2$, Co$_3$]: a ring $\,T\,$ is {\it projective-free\/} if
every finitely generated projective right module over $\,T\,$ is
free of a unique rank. (This notion is known to be left-right
symmetric; see [Co$_2$].) Clearly, such a ring $\,T\,$ has only
trivial idempotents (so it is abelian), and any idempotent in
$\,{\mathbb M}_n(T)\,$ is similar to a diagonal matrix with
$\,0$'s and $\,1$'s on the diagonal. 

\bs\nt
{\bf (3.10)} Another example of a ring satisfying the hypothesis
(2) in Theorem 3.8 is an abelian (von Neumann) regular ring
$\,T$. The fact that any idempotent matrix over a regular
ring is diagonalizable can be shown in at least two ways:
first by using the refinement theorem (for finitely generated
projective modules) of Goodearl and Handelman in [GH:~(3.8)],
and second, by using [Go:~(2.6)] in conjunction with
[SG:~Theorem 9].

\bs\nt
{\bf (3.11)} In the formulation of ``Open Problem (47)'' in
Goodearl's book [Go:~p.\,349], there is an extensive list
of different classes of regular rings $\,T\,$ that are
classically known to be $\,{\mathbb M}_n$-unique for all
$\,n$. Thus, Theorem 3.8 is applicable to all such regular
rings, under the hypothesis (1).  The two best known
classes among those mentioned in [Go:~p.\,349] are right
self-injective regular rings [Go:~(10.35)], and all regular
rings whose primitive factors are artinian [Go:~(6.12)].
(The latter class includes all abelian regular rings, which
we have already mentioned in (3.10).)  But unfortunately,
unit-regular rings are {\it not\/} known to be
$\,{\mathbb M}_n$-unique (even for $\,n=2$), so the answer
to Question 3.7 has so far remained unknown in the case
where the base ring $\,T\,$ is unit-regular.

\bs\nt
{\bf (3.12)} According to [La$_3$:~(17.26)], any semilocal
ring $\,T\,$ is $\,{\mathbb M}_n$-unique for all $\,n$. Thus,
Theorem 3.8 is applicable to $\,T\,$ again under the hypothesis
(1). In the special case where $\,(T, {\mathfrak m})\,$ is a
local ring, it would be projective-free as well. In this case,
as we have noted in the statement of Theorem B in \S1, the
conclusion $\,1_T\in {\rm U}(T)+{\rm U}(T)\,$ in Theorem 3.8
would amount to the simpler statement that
$\,|\,T/{\mathfrak m}\,|>2$. 

\bs\nt
{\bf (3.13)} Theorem 3.8 is applicable to any commutative ring
$\,T\,$ (under hypothesis (1)) since $\,T\,$ is also known to
be $\,{\mathbb M}_n$-unique (for all $\,n$) by [La$_3$:~(17.31)].
In this case, the conclusion of Theorem 3.8 is capable of
another more concrete derivation, as follows.  If $\,e,e'\,$
are idempotents in $\,R={\mathbb M}_2(T)\,$ with $\,\minus
\in {\rm U}(R)$, it is easy to show using [KLS:~Formula (1.1)]
that $\,{\rm tr}\,(ee')\,$ and $\,1-{\rm tr}\,(ee')\,$ are
both units in $\,T$, with sum $\,1$.

\bs
Using Theorem 3.3 and Corollary 3.5, it is now easy to
determine all matrix rings over local rings and nonzero
commutative rings that have Property K.

\bs\nt
{\bf Theorem 3.14.} {\it Let $\,T\,$ be a local ring or
a nonzero commutative ring. Then $\,R\!:={\mathbb M}_n(T)\,$
has {\rm Property K\/} iff either $\,n\in \{4,6,8,\dots\,\}$, 
or $\,n=2\,$ and $\,1_T\in {\rm U}(T)+{\rm U}(T)$.} 

\bs\nt 
{\bf Proof.} The ``if'' part follows from Theorem 3.3 and
Corollary 3.5.  For the ``only if'' part, assume that
$\,R={\mathbb M}_n(T)\,$ has Property K. Since $\,T\,$ is a
local ring or a nonzero commutative ring, there exists an ideal
$\,{\mathfrak m}\subseteq T\,$ such that $\,T/{\mathfrak m}\,$
is a division ring. Then $\,{\mathbb M}_n\bigl(T/{\mathfrak m}
\bigr)\cong {\mathbb M}_n(T)/{\mathbb M}_n({\mathfrak m})\,$
also has Property K. Applying Theorem 3.3, we have
$\,{\mathbb M}_n\bigl(T/{\mathfrak m}\bigr)\cong
{\mathbb M}_2(S)\,$ for some ring $\,S\,$ with
$\,1_S \in {\rm U}(S)+{\rm U}(S)$. By the uniqueness part
of Wedderburn's theorem, we see that $\,n\,$ must be even.
Finally, in the special case $\,n=2$, $\,R={\mathbb M}_2(T)\,$
having Property K would imply that $\,1_T\in {\rm U}(T)
+{\rm U}(T)\,$ by the remarks in (3.12) and (3.13).  \qed

\bs\nt
{\bf Example 3.15.} Let $\,R={\mathbb M}_n(D)\,$ be a typical
simple artinian ring, where $\,D\,$ is a division ring.
According to Theorem 3.14, $\,R\,$ has Property K except
precisely when $\,n\,$ is odd or when $\,n=2\,$ and
$\,|D|=2$.

\bs\nt
{\bf Example 3.16.} Let $\,{\mathbb Z}_{(p)}\,$ denote the
localization of $\,{\mathbb Z}\,$ at the prime ideal $\,(p)$.
If $\,T={\mathbb Z}\,$ or $\,{\mathbb Z}_{(2)}$, $\,{\mathbb M}_n
(T)\,$ has Property K iff $\,n\in \{4,6,8,\dots\,\}$. On the other
hand, if $\,T={\mathbb Z}_{(p)}\,$ where $\,p\,$ is an odd prime,
then $\,{\mathbb M}_n(T)\,$ has Property K iff $\,n\,$ is even.

\bs\nt
{\bf Example 3.17.} If the ring $\,T\,$ in Theorem 3.14 is
allowed to be noncommutative and not a local ring, the ``only
if" part in the theorem may no longer hold. For instance, we
may take $\,T\,$ to be any nonzero ring with Property K (so
$\,T\,$ is necessarily noncommutative). Then for {\it any}
$\,n\geq 1\,$ (odd or even), $\,R={\mathbb M}_n(T)\,$ also
has Property K since $\,T\,$ does.

\bs\nt
{\bf \S4. \ Interplay Between Property K and
Property $\boldsymbol{\overline{{\rm K}}}$}

\bs
In this section, we turn our attention to rings with Property
$\overline{{\rm K}}$. Our goal is to understand more precisely the
relationship and interaction between Property K and Property
$\overline{{\rm K}}$. We begin by recalling that any ring $\,R\,$
with $\,2\in {\rm U}(R)\,$ has Property $\overline{{\rm K}}$, but
not necessarily Property K. On the other hand, a $\,2\times 2\,$
matrix ring such as $\,R={\mathbb M}_2\bigl({\mathbb F}_2\bigr)\,$
does not have Property K (e.g.~by Theorem 3.3), and hence also does
not have Property $\overline{{\rm K}}$ since $\,{\rm char}\,(R)=2$.
From this example, it follows for instance that if $\,T\,$ is any
ring with an ideal of index $2$, then $\,{\mathbb M}_2(T)\,$ (with
factor ring $\,R\,$) does not have Property $\overline{{\rm K}}$,
and hence also not Property K.

\bs 
Our first main result in this section is Theorem 4.1 below, where
part (1) gives some necessary conditions on the invertibility
of an anti-commutator $\,\plus \in R$, while part (2) offers a
somewhat unexpected characterization for the rings with Property
$\overline{{\rm K}}$ without the use of products of idempotents.

\bs\nt
{\bf Theorem 4.1.} 
(1) {\it If $\,e,e'\in {\rm idem}\,(R)\,$ are such that $\plus 
\in {\rm U}(R)$, then $\,e'\,$ is similar to $\,e$, and we have 
$\,\langle \,e,w\,\rangle \in {\rm U}(R)\,$ for some $\,w\in 
{\rm U}(R)$.}\\
(2) {\it A ring $\,R\,$ has {\rm Property $\overline{{\rm K}}$\/}
iff there exists a unit $\,u\in {\rm idem}\,(R)+{\rm idem}\,(R)\,$
such that $\,1-u \in {\rm U}(R)$.}
 
\bs\nt
{\bf Proof.} (1) According to Theorem 2.13(2), $\plus \in
{\rm U}(R)\,$ amounts to two conditions: 
$$
e+e'\in {\rm U}(R), \;\;\mbox{and}\;\; 1-e-e'\in {\rm U}(R).
\leqno (4.2)
$$ 
The latter implies (by Theorem 2.8) that $\,e'=wew^{-1}\,$ 
for some $\,w \in {\rm U}(R)$.  Then $\,e+wew^{-1}=e+e'\in
{\rm U}(R)\,$ implies that $\,\langle\,e,w \,\rangle=ew+we=uw
\in {\rm U}(R)$.

\mk\nt
(2) If $\,R\,$ has Property $\overline{{\rm K}}$, there exist
$\,e,e'\in {\rm idem}\,(R)\,$ satisfying (4.2). Adding the two
elements in (4.2) shows that $\,1\,$ is a sum of two units, the
first one of which is a sum of two idempotents in $\,R$. Conversely,
if $\,1=u+u'\,$ where $\,u,u'\in {\rm U}(R)\,$ and $\,u=e+e'\,$
for some $\,e,e' \in {\rm idem}\,(R)$, then $\,1-e-e'=1-u=u'
\in {\rm U}(R)$. Thus, $(4.2)$ is satisfied, so $\,R\,$ has
Property $\overline{{\rm K}}$.\qed

\bs
Theorem 4.1 for rings with Property $\overline{{\rm K}}$ has a
complete analogue for rings with Property K too, which we shall
present below. One good reason we have chosen to prove Theorem
4.1 for Property $\overline{{\rm K}}$ first is that the proof
of its part (2) is easier and more intuitive, while it gives
an impetus toward finding an analogous (but somewhat harder)
characterization result in the case of Property K.

\bs\nt
{\bf Theorem 4.3.}
(1) {\it If $\,e,e'\in {\rm idem}\,(R)\,$ are such that $\,\minus
\in {\rm U}(R)$, then we have $\,[\,e,v\,] \in {\rm U}(R)\,$ for
some $\,v\in {\rm U}(R)$.}\\
(2) {\it A ring $\,R\,$ has {\rm Property K\/} iff there exists
a unit $\,v\in {\rm idem}\,(R)-{\rm idem}\,(R)\,$ such that
$\,1\,{\pm\,}v\in {\rm U}(R)\,$ for both signs.}
 
\bs\nt
{\bf Proof.} (1) By Theorem 2.13(1), $\,\minus \in {\rm U}(R)
\,\Rightarrow\, v\!:=e-e'\in {\rm U}(R)$, and so
$$
[\,e,v\,]=[\,e,\,e-e'\,]=-[\,e,e'\,]\in {\rm U}(R).  \leqno (4.4)
$$
The fact that one can ``replace'' the idempotent $\,e'\,$ by
a unit $\,v\,$ (while still achieving the property $\,[\,e,v\,]
\in {\rm U}(R)$) reminds us somewhat of the theorem of de S\'a
as discussed in [Sa] and [KL], although our result here applies
quite generally to any ring $\,R$, as long as $\,e,e' \in
{\rm idem}\,(R)$. (Needless to say, the converse of (1) fails
in general. For instance, taking $\,R={\mathbb M}_2(T)\,$ for
any ring $\,T$, $\,[e,v]\,$ is an invertible commutator for
$\,e=E_{11}\in {\rm idem}\,(R)\,$ and $\,v\in {\rm U}(R)$. But
for $\,S={\mathbb Z}\,$ or $\,{\mathbb F}_2\,$ for instance,
we cannot find {\it any\/} invertible $\,[e_1,e_2]\,$ with
$\,e_1,e_2\in {\rm idem}\,(R)$.)

\mk\nt
(2) First assume $\,R\,$ has Property K and fix an invertible
commutator $\,\minus\,$ with $\,e,e'\in {\rm idem}\,(R)$. We
have $\,(1-e)-e'\in {\rm U}(R)\,$ (by (2.13)(1)), so Theorem
2.8 implies that $\,(1-e)+e'\in {\rm U}(R)$; that is, $\,1-v\in
{\rm U}(R)\,$ for $\,v\!:=e-e'\in {\rm U}(R)\,$ as in the last
paragraph. On the other hand, $\,1+v=(1-e')+e\in {\rm U}(R)\,$
too, again by Theorem 2.8 since $\,(1-e')-e\in {\rm U}(R)$.
This proves the ``only if'' part of (2). For the converse,
assume that there is a unit $\,v=e-e'\,$ for suitable $\,e,e'
\in {\rm idem}\,(R)\,$ such that $\,1 {\pm\,}v \in {\rm U}(R)\,$
for both signs. By Theorem 2.13(1), we will have Property K on
$\,R\,$ if we can show that $\,1-e-e'\in {\rm U}(R)$. Appealing
to Theorem 2.8 once more (with the idempotent $\,e\,$ there
replaced by $\,1-e$), this would follow if $\,(1-e)+e'\in
{\rm U}(R)\,$ and $\,e+(1-e')\in {\rm U}(R)$. The former
amounts to $\,1-v\in {\rm U}(R)$, while the latter amounts
to $\,1+v\in {\rm U}(R)$. Since both of these conditions were
given, we are done.\qed

\bs
Our next goal is to study the possible presence of Property
$\overline{{\rm K}}$ on $\,n\times n\,$ matrix rings over
certain types of rings $\,T$. The first case we can treat
without too much difficulty is when $\,(T, {\mathfrak m})\,$
is a local ring. Here, a judicious use of Corollary 2.16
enables us to (essentially) ``replace'' $\,T\,$ by the
division ring $\,T/{\mathfrak m}$. As we have pointed out
before in (3.9), the condition $\,1_T\in {\rm U}(T)
+{\rm U}(T)\,$ in the statement (2) below can be more simply
expressed by $\,|\,T/{\mathfrak m}\,|>2$. However, we still
prefer to use the former condition because of its more general
nature. Similarly, in the statement (3) below, we prefer the
condition ``$2\in {\rm U}(T)$'' to the equivalent condition
``${\rm char}\,\bigl(T/{\mathfrak m}\bigr)\neq 2\,$''. 

\bs\nt
{\bf Theorem 4.5.} {\it Let $\,(T,{\mathfrak m})\,$ 
be a local ring. Then $\,R={\mathbb M}_n(T)\,$ has
{\rm Property $\overline{{\rm K}}$\/} iff we are in
one of the following cases\/}:

\mk\nt
(1) $\,n\in \{4,6,8,\dots \,\}$.\\
(2) $\,n=2\,$ {\it and $\,1_T\in {\rm U}(T)+{\rm U}(T)$.} \\
(3) $\,n\,$ {\it is odd and $\,2\in {\rm U}(T)$.}

\bs\nt
{\bf Proof.} If (1) or (2) holds, we know from Corollary 3.5
and Theorem 3.3 respectively that $\,R\,$ has Property K,
so of course $\,R\,$ has Property $\overline{{\rm K}}$. 
If (3) holds instead, then $\,<\!\!I_n,I_n\!\!>=2\,I_n
\in {\rm U}(R)\,$ shows that $\,R\,$ has Property
$\overline{{\rm K}}$. {\it Conversely, assume in the
following that $\,R\,$ has\/} Property $\overline{{\rm K}}$.
If $\,n\in \{4,6,8,\dots \,\}$, then (1) holds, so for the rest
of the proof, we need only work with odd $\,n\,$ and $\,n=2$.

\mk\nt
{\bf Case A.} {\it $\,T\,$ is a division ring.} Fix an
anti-commutator $\,\plus\in {\rm U}(R)\,$ with (necessarily
nonzero) $\,e,e'\in {\rm idem}\,(R)$. {\it First assume that
$\,n=2$.}  If $\,e=I_2$, $\,\plus =2\,e'\in {\rm U}(R)\,$
implies that $\,2\in {\rm U}(T)$, in which case $\,1_T
=2^{-1}+2^{-1}\in {\rm U}(T)+{\rm U}(T)$. If $\,e\neq I_2$,
we may assume (after a conjugation) that $\,e'=E_{11}$. Writing
$\,e=\footnotesize{\begin{pmatrix}a&b\\c&d\end{pmatrix}}$,
we have $\,\plus =\footnotesize{\begin{pmatrix}2a&b\\c&0
\end{pmatrix}}\in {\rm U}(R)$, so $\,b,c\in {\rm U}(T)$.
As $\,e^2=e\,$ implies that $\,a^2+b\,c=a$, we see that
$\,a,\,1-a\,$ are in $\,{\rm U}(T)\,$ too, with sum $\,1$,
so (2) holds. {\it Finally, suppose $\,n\,$ is odd.}  Here we
may assume that $\,e'={\rm diag}\,(I_k,0_{n-k})\,$ where $\,1
\leq k \leq n$. If $\,2\notin {\rm U}(T)$, then $\,2=0\in T$,
and so $\,k<n\,$ (for otherwise $\plus =2\,e'=0$, which is
impossible).  Writing $\,e=\footnotesize{\begin{pmatrix}A&B\\
C&D \end{pmatrix}}$ with $\,A\in {\mathbb M}_k(T)\,$ and
$\,D \in {\mathbb M}_{n-k}(T)$, we have
$\,\plus=\footnotesize{\begin{pmatrix}0&B\\C&0\end{pmatrix}}$.
As $\,B\,$ and $\,C\,$ are {\it non-square matrices}, we can take
a nonzero vector $\,v\,$ such that either $\,B\,v=0\,$ or $\,C\,v=0$.
Then we'll have either $\,\plus \footnotesize{\begin{pmatrix}
0\\v\end{pmatrix}}=0\,$ or $\,\plus \footnotesize{\begin{pmatrix}
v\\0\end{pmatrix}}=0$, in contradiction to the fact that $\,\plus
\in {\rm U}(R)$. This completes the proof that $\,2\in
{\rm U}(T)\,$ when $\,n\,$ is odd.

\mk\nt
{\bf Case B.} {\it $(T,{\mathfrak m})\,$ is a local ring.}
Let $\,\overline{T}\,$ be the division ring
$\,T/{\mathfrak m}$.  By [La$_1$:~p.\,57], 
$\,{\rm rad}\,({\mathbb M}_n(T))={\mathbb M}_n({\rm rad}\,(T))
= {\mathbb M}_n({\mathfrak m})$. Therefore,
$$
{\mathbb M}_n(T)/{\rm rad}\,({\mathbb M}_n(T))
={\mathbb M}_n(T)/{\mathbb M}_n({\mathfrak m})
\cong {\mathbb M}_n\bigl(\overline{T}\bigr). \leqno (4.6)
$$
Since $\,{\mathbb M}_n(T)\,$ is an exchange ring (by [Wa] or [Ni]),
Corollary 2.16 implies that it has Property $\overline{{\rm K}}$
iff $\,{\mathbb M}_n\bigl(\overline{T}\bigr)\,$ does. Therefore,
we are free to replace $\,T\,$ by $\,\overline{T}\,$ to assume
that $\,T\,$ is a division ring, in which case we are fully
covered by Case A above. \qed

\bs
In general, we do not know exactly when a general matrix ring
$\,{\mathbb M}_n(T)\,$ will have Property $\overline{{\rm K}}$
if $\,n\,$ is odd or $\,n=2$. Aside from the case where $\,T\,$
is local (as treated above in Theorem 4.5), a second manageable
case is where $\,T\,$ is a commutative ring. Working under this
assumption, it turns out that the more substantial case to
consider is $\,n=2$. In this case, we were pleasantly surprised
to find that Property K and Property $\overline{{\rm K}}$ are
{\it equivalent\/} on $\,{\mathbb M}_2(T)\,$!  To prove this,
we start by first working out some easy facts about determinants
of $\,2\times 2\,$ matrices over a commutative ring.

\bs\nt
{\bf Proposition 4.7.} {\it For $\,A,B\in R={\mathbb M}_2(T)\,$
where $\,T\,$ is a commutative ring, we have}
$$
{\rm det}\,(A+B)+{\rm det}\,(A-B)
=2\,[\,{\rm det}\,(A)+{\rm det}\,(B)\,] \in T. \leqno (4.8)
$$

\nt
{\bf Proof.} Writing $\,A=\footnotesize{\begin{pmatrix}p&q\\r&s
\end{pmatrix}}$ and $\,B=\footnotesize{\begin{pmatrix}w&x\\y&z
\end{pmatrix}}$, a quick computation shows that
$$
{\rm det}\,(A+B)={\rm det}\,(A)+{\rm det}\,(B)+(pz+sw-qy-rx),
\;\mbox{and}
$$
\vspace{-7.5 mm}
$$
{\rm det}\,(A-B)={\rm det}\,(A)+{\rm det}\,(B)-(pz+sw-qy-rx).
$$
Adding these two formulas gives the desired equation (4.8).\qed

\bs\nt
{\bf Corollary 4.9.} {\it In the notations of \,{\rm Proposition
4.7}, if $\,2\,[\,{\rm det}\,(A)+{\rm det}\,(B)\,] \in
{\rm rad}\,(T)$, then $\,A+B\in {\rm U}(R)\,$ iff $\,A-B \in
{\rm U}(R)$. In particular, this} ``{\it iff\/}'' {\it statement
holds in case $\,{\rm det}\,(A)+{\rm det}\,(B)=0$; or more
specifically, if} $\,{\rm det}\,(A)={\rm det}\,(B)=0$.
  
\bs\nt
{\bf Proof.} The first ``iff'' statement holds on account
of (4.8) since a matrix $\,C\in R\,$ is invertible iff
$\,{\rm det}\,(C)\in {\rm U}(T)$, while for any $\,t,\,t'\in T$,
$\,t+t' \in {\rm rad}\,(T)\,$ implies that $\,t\in {\rm U}(T)\,$
iff $\,t'\in {\rm U}(T)$. The rest of the Corollary is clear.\qed

\bs
With the above Corollary providing a crucial link between the
invertibility of $\,A+B\,$ and $\,A-B$, we are now ready to
prove the following result. 

\bs\nt
{\bf Theorem 4.10.} {\it Let $\,R={\mathbb M}_2(T)\,$ where
$\,T\,$ is a commutative ring. Then $\,R\,$ has {\rm Property K}
iff it has {\rm Property\/} $\overline{{\rm K}}$.}

\bs\nt
{\bf Proof.} Of course, only the ``if'' part is at stake. Before
beginning its proof, we first point out that the ``if'' part
does not mean that for any $\,e,e'\in R$, $\,\plus \in
{\rm U}(R)\,$ implies that $\,\minus \in {\rm U}(R)$.  Indeed,
in the case where $\,2\in {\rm U}(T)$, choosing $\,e=e'=I_2\,$
gives $\,\plus =2\,I_2 \in {\rm U}(R)$, but $\,\minus =0
\notin {\rm U}(R)\,$ if $\,T\neq 0$. (Nevertheless, under
the assumption that $\,2\in {\rm U}(T)$, $\,R\,$ does have
Property K by taking $\,a=b=2^{-1}\,$ in the proof of
Theorem 3.1.) {\it In the following, we assume that $\,R\,$
has Property $\overline{{\rm K}}$.} We'll show that $\,R\,$
has Property K in two steps.

\mk\nt
{\bf Step 1.} {\it The desired conclusion is true if $\,T\,$ is
a connected ring\/}; that is, if $\,{\rm idem}\,(T)=\{0,1\}$. 
Indeed, taking an invertible $\,\plus\,$ (for suitable $\,e,e'
\in {\rm idem}\,(R)$), $\,{\rm det}\,(e)\,$ is either $\,1\,$
or $\,0$. In the former case, $\,e\,$ is invertible, so
$\,e=I_2$. Here, $\,\plus=2\,e'\in {\rm U}(R)\,$ implies that
$\,2\in {\rm U}(T)$, so of course $\,R\,$ has Property K.
We may thus assume that $\,{\rm det}\,(e)=0$, in which case
$\,{\rm det}\,(ee')={\rm det}\,(e'e)=0$. As $\,ee'+e'e\in
{\rm U}(R)$, the last part of Corollary 4.9 implies that
$\,ee'-e'e\in {\rm U}(R)$, so $\,R\,$ has Property K.

\mk\nt
{\bf Step 2.} {\it Assume that $\,R\,$ does not have} Property K.
By applying Zorn's Lemma to the family $\,{\mathcal F}\,$ of
ideals $\,J_i\subseteq T\,$ such that $\,{\mathbb M}_2(T/J_i)\,$
does not have Property K, we see that $\,{\mathcal F}\,$ has
a maximal member (with respect to inclusion), say $\,J$.  The
maximal choice of $\,J\,$ implies that $\,T/J\,$ is a connected
ring, for otherwise $\,{\mathbb M}_2(T/J)\,$ would have been a
direct product of a pair of $\,2\times 2\,$ matrix rings each
of which has Property K, in contradiction to the fact that $\,J
\in {\mathcal F}$.  On the other hand, $\,{\mathbb M}_2(T/J)\,$
is isomorphic to a factor ring of $\,{\mathbb M}_2(T)$, so
it has Property $\overline{{\rm K}}$. As $\,T/J\,$ is a
connected ring, what we have done in Step 1 implies that
$\,{\mathbb M}_2(T/J)\,$ has Property K. This is a
contradiction. \qed

\bs
With the help of Theorem 4.10, we can now decide exactly when 
a matrix ring $\,{\mathbb M}_n(T)\,$ over a commutative ring
$\,T\,$ will have Property $\overline{{\rm K}}$.  It is of
interest to note that the conclusions in the following theorem
happen to be {\it identical\/} to those of Theorem 4.5 over
a (possibly noncommutative) local ring $\,T$, although the
proofs are rather different in the key case where $\,n=2$. 

\bs\nt
{\bf Theorem 4.11.} {\it For any commutative ring $\,T$, $\,R
={\mathbb M}_n(T)\,$ has {\rm Property $\overline{{\rm K}}$\/}
iff we are in one of the following cases\/}:

\mk\nt
(1) $\,n\in \{4,6,8,\dots \,\}$.\\
(2) $\,n=2\,$ {\it and $\,1_T\in {\rm U}(T)+{\rm U}(T)$.} \\
(3) $\,n\,$ {\it is odd and $\,2\in {\rm U}(T)$.}

\bs\nt
{\bf Proof.} If (1) or (2) holds, $\,R\,$ would have Property K,
and therefore Property $\overline{{\rm K}}$. If (3) holds, of
course $\,R\,$ has trivially Property $\overline{{\rm K}}$.
{\it Conversely, assume $\,R\,$ has} Property $\overline{{\rm K}}$.
If $\,n\in \{4,6,8,\dots \,\}$, we are in Case (1). If $\,n=2$,
we know from Theorem 4.10 that $\,R\,$ has Property K, so
Theorem 3.14 implies that $\,1_T\in {\rm U}(T)+{\rm U}(T)$.
Finally, let $\,n\,$ be odd. If $\,2\notin {\rm U}(T)$, then
$\,2\,$ lies in some maximal ideal $\,{\mathfrak m}\subseteq T$.
Since $\,{\mathbb M}_n\bigl(T/{\mathfrak m}\bigr)\,$ is isomorphic
to a factor ring of $\,R$, it also has Property $\overline{{\rm K}}$.
As $\,T/{\mathfrak m}\,$ is a field of characteristic $\,2$,
$\,{\mathbb M}_n\bigl(T/{\mathfrak m}\bigr)\,$ would also have
Property K, which would contradict Theorem 3.14. Therefore,
we must have $\,2\in {\rm U}(T)$, so we are in Case (3).\qed

\bs
The result above suggests that the case of $\,{\mathbb M}_2(T)\,$
for commutative rings $\,T\,$ is of special interest in the
treatment of Property K (or equivalently, Property
$\overline{{\rm K}}\,$ by Theorem 4.10). We'll now conclude this
paper by giving some more characterizations for these properties
by using commutators and anti-commutators whose second entries
are diagonalizable matrices instead of idempotent matrices. Such
a replacement is by no means automatic since in general these
two classes of matrices are logically independent.

\bs\nt
{\bf Theorem 4.12.} {\it For $\,R={\mathbb M}_2(T)\,$ where
$\,T\,$ is a commutative ring, the following statements
are equivalent\/}:

\mk\nt
(1) $\,R\,$ {\it has\/} Property K ({\it or equivalently,}
$\,1_T\in {\rm U}(T)+{\rm U}(T)$). \\
(2) $\,[\,e,\delta\,] \in {\rm U}(R)\,$ {\it for
some $\,e\in {\rm idem}\,(R)\,$ and some diagonalizable
matrix $\,\delta\in R$.}\\
(3) $\,\langle\,e,\delta\,\rangle \in {\rm U}(R)\,$ {\it for
some $\,e\in {\rm idem}\,(R)\,$ and some diagonalizable
matrix $\,\delta\in R$.}

\mk\nt
{\it In the case where $\,T\,$ is an exchange ring,
these statements are also equivalent to\/}:

\mk\nt
(4) $\,T\,$ {\it has no ideal of index $\,2$.}

\bs\nt
{\bf Proof.} First assume (1) holds; say $\,1=a+b\,$
for some $\,a,b\in {\rm U}(T)$. Taking
$\,e=\footnotesize{\begin{pmatrix}a&b\\a&b \end{pmatrix}} \in
{\rm idem}\,(R)\,$ and $\,\delta={\rm diag}\,(1,0)$, we have
$\,[\,e,\delta\,] \in \footnotesize{\begin{pmatrix}0&b\\-a&0
\end{pmatrix}}\in {\rm U}(R)$, and $\,\langle\,e,\delta\,\rangle
\in \footnotesize{\begin{pmatrix}2\,a&b\\a&0\end{pmatrix}}\in
{\rm U}(R)$, so (2) and (3) both hold.  Next, $\,(2)\Rightarrow
(1)$ follows from  the calculation in the proof of Case (2) in
Theorem 3.8, after a reduction to the case where $\,\delta\,$
is diagonal. Now we come to the harder implication
$\,(3)\Rightarrow (1)$.  Given $\,e,\,\delta\,$ as in (3), we
may again assume that $\,\delta={\rm diag}\,(s,t)\,$ for some
$\,s,t\in T$. Letting $\,e=\footnotesize{\begin{pmatrix}
a&b\\c&d\end{pmatrix}}$, the determinant of
$\,\langle\,e,\delta\,\rangle\,$ is easily computed to be
$$
4\,adst-b\,c\,(s+t)^2 \in {\rm U}(T). \leqno (4.13)
$$
Let $\,g={\rm det}\,(e)=ad-bc$, which is an idempotent in $\,T$.
We have projection maps from $\,T\,$ onto $\,gT\,$ and onto
$\,(1-g)\,T$.  Taken together, these define a natural ring
isomorphism from $\,T\,$ onto $\,gT \times (1-g)\,T$. Therefore,
our job is reduced to showing that $\,1\,$ is a sum of two units
in both $\,gT\,$ and $\,(1-g)\,T$. Under the two projections,
$\,g\,$ projects to the identity and to zero respectively.
Thus, it suffices to argue in the following two cases.

\mk\nt
{\bf Case 1.} $\,g=1$. In this case, $\,e=I_2$, so
$\,\langle\,e,\delta\,\rangle = 2\,\delta \in {\rm U}(R)\,$
implies that $\,2\in {\rm U}(R)$. Thus, $\,2\in {\rm U}(T)$,
and so $\,1=2^{-1}+2^{-1}\in {\rm U}(T)+{\rm U}(T)$.

\mk\nt
{\bf Case 2.} $\,g=0$. In this case, $\,ad=bc$, so (4.13) shows
that $\,b\,c\in {\rm U}(T)$.  But $\,e\in {\rm idem}\,(R)\,$
implies that $\,a^2+b\,c=a$, so $\,a,\,1-a\,$ are units in $\,T$,
with sum $\,1$.

\mk
Finally, $(1)\Rightarrow (4)\,$ is always true, since the
existence of an ideal $\,I\subseteq T\,$ of index two
would imply that $\,R\,$ has a factor ring isomorphic to
$\,{\mathbb M}_2\bigl({\mathbb F}_2\bigr)$, which does not
have Property K. Conversely, assume that $\,T\,$ is an
exchange ring satisfying (4). To prove (1), we assume 
instead that $\,1_T\notin {\rm U}(T)+{\rm U}(T)$. Consider
the nonempty family $\,{\mathcal G}\,$ of ideals $\,I_i
\subseteq T\,$ for which $\,1\,$ is not a sum of two units
in $\,T/I_i$.  Applying Zorn's Lemma to the family
$\,{\mathcal G}$, we see that $\,{\mathcal G}\,$ has a
maximal member $\,I\,$ (with respect to inclusion). The
maximal choice of $\,I\,$ implies (as in the proof of (4.10))
that the commutative factor ring $\,T/I\,$ is connected.
Since $\,T/I\,$ remains to be an exchange ring (and $\,T/I
\neq 0$), it must be a local ring (by [Wa:~Proposition 1]).
As $\,1\,$ is not a sum of two units in $\,T/I$, the local
ring $\,T/I\,$ must have residue field $\,{\mathbb F}_2$.
This implies that $\,T\,$ has a factor ring isomorphic to
$\,{\mathbb F}_2$, which contradicts (4). \qed

\mk

\bs\nt
{\bf \S5. \ Commutator Characterizations for
$\,{\boldsymbol 2}\times {\boldsymbol 2}\,$ Matrix Rings}

\bs
According to Theorem 3.3, rings with Property K are precisely
$\,2\times 2\,$ matrix rings over base rings with a specific
unit property. One may naturally ask: how about characterizations
of the {\it most general\/} $\,2\times 2\,$ matrix rings?
In the literature, there is a simple criterion for these
(${\rm (A)}\Leftrightarrow {\rm (H)}$ in Theorem 5.1), that
is due to Fuchs, Maxson and Pilz [FMP:~(III.2)]. To make our
exposition more useful and yet completely self-contained, we'll
reprove the Fuchs-Maxson-Pilz theorem in a much expanded form
below, following the idea of using anti-commutators (as well
as commutators) in [La$_4$:~p.\,349] but not assuming any
general matrix ring recognition theorems in [La$_3$:~\S17].
In the following result, the statements (B) through (G) are
expressed in terms of anti-commutators, while (I), (J) and (K)
are expressed in terms of commutators.

\bs\nt
{\bf Theorem 5.1.} {\it For any ring $\,R$, the following 
statements are equivalent\/}:

\sk\nt
(A) {\it $\,R\cong {\mathbb M}_2(S)\,$ for some ring $\,S$.}\\
(B) {\it There exist $\,p,q\in R\,$ with $\,p^2=q^2=0\,$ such
that $\,\langle\,p,q\,\rangle =1$.}\\
(C) {\it There exist $\,p,q\in R\,$ with $\,p^2=q^2=0\,$ such
that $\,\langle\,p,q\,\rangle \in {\rm U}(R)$.}\\
(D) {\it There exist $\,v\in {\rm U}(R)\,$ and $\,p\in R\,$
with $\,p^2=0\,$ such that $\,\langle\,p,v\,\rangle =1$.}\\
(E) {\it There exist $\,v\in {\rm U}(R)\,$ and $\,p\in R\,$
with $\,p^2=0\,$ such that $\,\langle\,p,v\,\rangle
\in {\rm U}(R)$.}\\
(F) {\it There exist $\,r,p\in R\,$ with $\,p^2=0\,$ such
that $\,\langle\,p,r\,\rangle =1$.}\\
(G) {\it There exist $\,r,p\in R\,$ with $\,p^2=0\,$ such
that $\,\langle\,p,r\,\rangle \in {\rm U}(R)$.}\\
(H) {\it There exist $\,p,q\in R\,$ with $\,p^2=q^2=0\,$ and
$\,p+q\in {\rm U}(R)$.}\\
(\,I\,) {\it There exist $\,p,q\in R\,$ with $\,p^2=q^2=0\,$
such that $\,[\,p,q\,]\in {\rm U}(R)$.}\\
(J\,) {\it There exist $\,v\in {\rm U}(R)\,$ and $\,p\in R\,$
with $\,p^2=0\,$ such that $\,[\,p,v\,]\in {\rm U}(R)$.}\\
(K) {\it There exist $\,r,p\in R\,$ with $\,p^2=0\,$ such
that $\,[\,p,r\,]\in {\rm U}(R)$.}

\bs\nt
{\bf Proof.} ${\rm (A)}\Rightarrow {\rm (B)}$. If (say)
$\,R={\mathbb M}_2(S)$, then (B) holds for $\,p=E_{12}\,$
and $\,q=E_{21}$.

\mk\nt
${\rm (B)}\Rightarrow {\rm (C)}$ is a tautology.

\mk\nt
${\rm (C)}\Rightarrow {\rm (H)}$. For $\,p,q\,$ as in (C),
we have $\,(p+q)^2=\langle\,p,q\,\rangle \in {\rm U}(R)$,
so $\,p+q\in {\rm U}(R)$. 

\mk\nt
${\rm (H)}\Rightarrow {\rm (D)}$. For $\,p,q\in R\,$ as in (H),
let $\,v\,$ be the inverse of $\,u\!:=p+q\in {\rm U}(R)$. Left
multiplying by $\,q\,$ and right multiplying by $\,p\,$ give
$\,qu=qp\,$ and $\,up=qp\,$ respectively. Equating these gives
$\,vq=pv$, so $\,\langle \,p,v\,\rangle =pv+vp=vq+vp=vu=1$.

\mk\nt
${\rm (D)}\Rightarrow {\rm (G)}\,$ is a tautology.

\mk\nt
${\rm (G)}\Rightarrow {\rm (F)}$. For $\,r,p\,$ as in (G),
let $\,w=pr+rp\in {\rm U}(R)$. Left multiplying by $\,p\,$
gives $\,pw=prp$, and right multiplying by $\,p\,$ gives
$\,wp=prp$. Thus, $\,pw=wp$, so left multiplying by
$\,w^{-1}\,$ yields $\,1=w^{-1}pr+w^{-1}rp=p\,(w^{-1}r)
+(w^{-1}r)\,p$, which verifies (F).

\mk\nt
${\rm (F)}\Rightarrow {\rm (A)}$. For $\,r,p\,$ as in (F),
left multiplying $\,pr+rp=1\,$ by $\,e\!:=rp\,$ shows that
$\,e\in {\rm idem}\,(R)$, with complementary idempotent
$\,f\!:=pr$. By [La$_1$:~(21.20)], $\,e,f\,$ are
{\it isomorphic\/} idempotents. Thus, we have $\,R_R=eR\oplus fR
\cong eR\oplus eR\,$ as right $\,R$-modules. Taking endomorphism
rings gives $\,R\cong {\rm End}_R(eR)\cong {\mathbb M}_2(eRe)$.

\mk\nt
${\rm (D)}\Leftrightarrow {\rm (E)}$. We need only show that
$\,{\rm (E)}\Rightarrow {\rm (D)}$, so let $\,w\!:=pv+vp
\in {\rm U}(R)\,$ where $\,p^2=0\,$ and $\,v\in {\rm U}(R)$.
Applying the same proof for ${\rm (G)}\Rightarrow {\rm (F)}\,$
here gives $\,1=p\,(w^{-1}v)+(w^{-1}v)\,p$. Since $\,w^{-1}v
\in {\rm U}(R)$, this proves (D).

\mk\nt
${\rm (A)}\Rightarrow \{\mbox{(I), (J) and (K)}\}$. If (say)
$\,R={\mathbb M}_2(S)$, then (I) holds for $\,p=E_{12}\,$ and
$\,q=E_{21}$, and (J), (K) both hold for $\,p=E_{12}\,$ and
$\,v=r=E_{12}+E_{21}$.

\mk\nt
${\rm (I)}\Rightarrow {\rm (C)}$. Given $\,p,q\,$ in (I), let
$\,w\!:=[\,p,q\,]\in {\rm U}(R)$. Repeating the first steps
in the proof of ${\rm (G)}\Rightarrow {\rm (F)}$, we can show
here that $\,w\,$ {\it anti-commutes\/} with both $\,p\,$
and $\,q$. Thus,
$$
1=w^{-1}pq-w^{-1}qp=-p\,(w^{-1}q)-(w^{-1}q)\,p
= -\langle\,p,w^{-1}q\,\rangle\,.   \leqno (5.2)
$$
Since $\,(w^{-1}q)^2=w^{-1}qw^{-1}q=-w^{-2}q^2=0$, this
verifies (C).

\mk\nt
${\rm (J)}\Rightarrow {\rm (K)}\Rightarrow {\rm (G)}$. It suffices
to prove the latter implication. Given $\,w\!:=pr-rp\in {\rm U}(R)\,$
as in (K), we have again $\,wp=-pw$. Thus, by repeating the
calculations in (5.2) (with $\,q\,$ replacing $\,r$), we get
$\,1=-\langle \,p,\,w^{-1}r\,\rangle\,$. Since $\,p^2=0$, this
verifies (G). \qed

\bs\nt
{\bf Remark 5.3.} It would be tempting to think that in the
statements (B), (D) and (F), we can replace the anti-commutators
by commutators and get three more equivalent statements in
terms of commutators.  However, this is not the case!
In fact, given (A), we may not be able to get an equation
$\,[\,p,q\,]=1\,$ with $\,p^2=q^2=0$. Indeed, if $\,R
={\mathbb M}_2(S)\,$ for a commutative ring $\,S$, {\it any\/}
commutator in $\,R\,$ will have zero trace, while the identity
matrix $\,I_2\,$ has trace $\,2$. Similarly, we must refrain
from trying to change the anti-commutator into a commutator
in the statements (D) and (F).  As for statement (H) (the
Fuchs-Maxson-Pilz characterization for (A)), it is easy to 
see that we also cannot change the requirement $\,p+q\in 
{\rm U}(R)\,$ into $\,p+q=1\in R\,$ (or $\,p+q\,$ being a 
{\it central\/} unit), since this would have implied that 
$\,p=q=0$, and hence $\,R=0\,$!

\bs\nt
{\bf Remark 5.4.} As we have pointed out in the paragraph
preceding Theorem 5.1, the proof of that theorem was designed
to be entirely self-contained so that it can be read
independently of the standard recognition theorems for
matrix rings as developed, for instance, in [La$_3$].
In the simplest case, according to [La$_3$:~(17.10)],
a standard characterization for a ring $\,R\,$ to be a
$\,2\times 2\,$ matrix ring is that there exist $\,p,r,r'
\in R\,$ with $\,p^2=0\,$ such that $\,pr+r'p=1$.  This
condition can thus be added to all the other equivalent
conditions listed in Theorem 5.1.  We note that it is
ostensibly ``weaker'' (as a sufficient condition
for (A)) than the three conditions (B), (D), and (F).

\bs
At this point, it would be helpful to recall a key definition
introduced in our earlier work [KL] on invertible commutators.
In that paper, we defined a ring element $\,a\in R\,$ to be
{\it completable\/} if there exists $\,r\in R\,$ such that
$\,[\,a,r\,] \in {\rm U}(R)$. Using this terminology, the
equivalence of (A) and (K) in Theorem 5.1 says precisely that
$\,R\,$ is a $\,2\times 2\,$ matrix ring (over some other ring)
iff $\,R\,$ has a square-zero completable element. Inspired by
this fact (as well as by our earlier Theorem 3.3), we stumbled
upon the idea of looking for {\it completable idempotents\/}
in $\,R$. It was a pleasant surprise to us that the existence
of a completable idempotent in a ring $\,R\,$ turns out to be
{\it also\/} equivalent to $\,R\,$ being a $\,2\times 2\,$
matrix ring. This fact and some of its variations are
collected in the new characterization theorem below for
$\,2\times 2\,$ matrix rings. Recall from Corollary 3.4 that,
in the special case where $\,2\in {\rm U}(R)$, we did know
that $\,R\,$ is a $\,2\times 2\,$ matrix ring iff $\,R\,$ has
Property K (that is, $\,\minus \in {\rm U}(R)\,$ for some
$\,e,e'\in {\rm idem}\,(R)$). However, the following new
equivalence result, like Theorem 5.1, is valid for all rings. 

\mk\nt
{\bf Theorem 5.5.} {\it For any ring $\,R$, the following 
statements are equivalent\/}:

\sk\nt
(1) {\it $\,R\cong {\mathbb M}_2(S)\,$ for some ring $\,S$.}\\
(2) {\it There exist $\,e\in {\rm idem}\,(R)\,$ and $\,r\in 
R\,$ with $\,r^2=1\,$ such that $\,[\,e,r\,]
\in {\rm U}(R)$.}\\
(3) {\it There exist $\,e\in {\rm idem}\,(R)\,$ and $\,r\in 
R\,$ with $\,r^2=-1\,$ such that $\,[\,e,r\,]\in {\rm U}(R)$.}\\
(4) {\it There exist $\,e\in {\rm idem}\,(R)\,$ and $\,r\in 
R\,$ with $\,r^3=1\,$ such that $\,[\,e,r\,]
\in {\rm U}(R)$.}\\
(5) {\it There exist $\,e\in {\rm idem}\,(R)\,$ and 
$\,r\in {\rm U}(R)\,$ such that $\,[\,e,r\,]\in {\rm U}(R)$.}\\
(6) {\it $\,R\,$ has a completable idempotent.}\\
(7) {\it There exist $\,e\in {\rm idem}\,(R)\,$ and $\,r\in R\,$ 
with $\,r^2=0\,$ such that $\,[\,e,r\,]\in {\rm U}(R)$.}\\
(8) {\it There exist $\,e\in {\rm idem}\,(R)\,$ and 
a nilpotent $\,r\in R\,$ such that $\,[\,e,r\,]\in {\rm U}(R)$.}
 
\mk\nt
{\bf Proof.} We first prove that (1) implies (2), (3), (4)
and (7). Indeed, if $\,R={\mathbb M}_2(S)\,$ for some ring $\,S$,
we can take $\,e=\footnotesize{\begin{pmatrix}1&0\\0&0\end{pmatrix}}$,
$\,r_2=\footnotesize{\begin{pmatrix}0&1\\1&0\end{pmatrix}}$,
$\,r_3=\footnotesize{\begin{pmatrix}0&1\\-1&0\end{pmatrix}}$, 
$\,r_4=\footnotesize{\begin{pmatrix}-1&1\\-1&0\end{pmatrix}}$, 
and $\,r_7=\footnotesize{\begin{pmatrix}1&1\\-1&-1\end{pmatrix}}$.
Then $\,e^2=e$, $\,r_2^2=r_4^3=1$, $\,r_3^2=-1$, and $\,r_7^2=0$,
and we check easily that $\,[\,e,r_i\,]\in {\rm U}(R)\,$ for all $\,i$. 

\mk\nt
$(7)\Rightarrow (8) \Rightarrow (6)\,$ are tautologies. [We note
on the side that (7) also trivially implies (K) (and hence (A))
in Theorem 5.1, although this information is not needed here.]

\mk\nt
$(6)\Rightarrow (1)$. Let $\,e=e^2\in R\,$ and $\,r\in R\,$
be such that the following commutator is a unit:
$$
[\,e,r\,]=er-re=er\,(1-e)-(1-e)\,re. \leqno (5.6)
$$
Letting $\,p=er\,(1-e)\,$ and $\,q=-(1-e)\,re$, we have 
clearly $\,p^2=q^2=0$, and $\,p+q\in {\rm U}(R)\,$ according
to (5.6). Thus, $\,{\rm (H)}\Rightarrow {\rm (A)}\,$ in
Theorem 5.1 gives (1).

\mk
The proof is now complete, since $\{(2)\;\mbox{or} \;(3)\;\mbox{or}
\;(4)\} \Rightarrow (5)\Rightarrow (6)\,$ are tautologies.\qed

\bs
It is of interest to point out that, in the various criteria above,
$\,r\,$ can be chosen to be a unit (respectively, a square root
of $\,{\pm 1}\,$ or a cubic root of $\,1$), or a nilpotent element
(respectively, a square-zero element), {\it but not necessarily an
idempotent.} We note also that, in the criterion (5) above, the
commutator $\,[\,e,r\,]\,$ cannot be replaced by the anti-commutator
$\,\langle\,e,r\,\rangle$.  Indeed, if we take $\,R\,$ to be
any nonzero commutative ring with $\,2\in {\rm U}(R)$, and
choose $\,e=r=1$, then $\,\langle\,e,r\,\rangle \in {\rm U}(R)$.
Here, $\,R\,$ cannot possibly be a $\,2\times 2\,$ matrix ring.
Nevertheless, we do have the following result for anti-commutators.

\bs\nt
{\bf Corollary 5.7.} {\it If $\,R\,$ is a ring with an idempotent
$\,e\,$ and an element $\,r\,$ such that $\,\langle \,e,r\,\rangle
\in {\rm U}(R)$, then $\,R/2R\cong {\mathbb M}_2(S)\,$ for some
ring $\,S$.}

\bs\nt
{\bf Proof.} In $\,\overline{R}=R/2R$, we have $\,\overline{er-re}
=\overline{er+re}\in {\rm U}\bigl(\overline{R}\bigr)$, so we can
apply $\,(5)\Rightarrow (1)\,$ in Theorem 5.5 to the factor ring 
$\,\overline{R}$.\qed

\bs
Emboldened by the fact that $\,2\times 2\,$ matrix rings are
characterized by the existence of completable idempotents as
well as the existence of completable square-zero elements, we
were led to the consideration of {\it completable involutions\/}
as well.  (By an involution, we simply mean an element $\,u\in
R\,$ with $\,u^2=1$.) In ring theory, it is rather rare that
idempotents, involutions and square-zero elements would play
parallel roles in the treatment of a certain problem or property.
But for the problem of characterizing $\,2\times 2\,$ matrix
rings, this does turn out to be the case, as the following
result shows.

\bs\nt
{\bf Theorem 5.8.} {\it A ring $\,R\,$ is a $\,2\times 2\,$ matrix
ring iff it has a completable involution.} 

\bs\nt
{\bf Proof.} If $\,R={\mathbb M}_2(S)\,$ for some ring $\,S$, we
have already pointed out earlier that the involution $\,E_{12}
+E_{21}\in R\,$ is completable. Conversely, assume that there exist  
$\,u,r\in R\,$ with $\,u^2=1\,$ such that $\,v\!:=[\,u,r\,]\in
{\rm U}(R)$. Repeating the idea (in the proof of (5.1)) of left and
right multiplying by $\,u$, and keeping in mind the equation
$\,u^2=1$, we see here that $\,uv=-vu$. Thus,
$$
1=v^{-1}ur-v^{-1}ru=(1-u)\,v^{-1}r-v^{-1}r\,(1+u); \leqno (5.9)
$$
that is, $\,1=e+f$, where $\,e\!:=(1-u)\,v^{-1}r\,$ and
$\,f\!:=-v^{-1}r\,(1+u)$. Since $\,fe=0$, we have $\,f\in
{\rm idem}\,(R)$, with complementary idempotent $\,e$. {\it It
suffices to show that $\,e\,$ and $\,f\,$ are isomorphic idempotents,}
as that will show (as in $\,{\rm (F)}\Rightarrow {\rm (A)}\,$ in
the proof of Theorem 5.1) that $\,R\cong {\mathbb M}_2(eRe)$. Right
multiplying (5.9) by $\,1-u\,$ and left multiplying it by $\,1+u\,$
give the following two ``von Neumann regularity'' equations:
$$
(1-u)\,v^{-1}r\,(1-u)=1-u,\;\;\mbox{and}\;\;\,
(1+u)\,(-v^{-1}r)\,(1+u)=1+u. 
$$
From the former, we have $\,eR=(1-u)\,R$. From the latter, we
have $Rf=R\,(1+u)$, which is well known to imply that $\,fR\cong
(1+u)\,R\,$ as right $\,R$-modules. (For a full proof of this
implication, see [La$_2$:~Exercise 1.17].)  Noting that
$$
(1-u)\,R=(1-u)\,vR=v\,(1+u)\,R \cong (1+u)\,R,
$$
we conclude that $\,eR\cong fR$, as desired.\qed

\bs\nt
{\bf Remark 5.10.} In view Theorem 5.8 (and parts of Theorem 5.5),
it might be tempting to surmise that the existence of a completable
element $\,u\in R\,$ with $\,u^n=1\,$ for some $\,n\geq 3\,$ might
also imply that $\,R\,$ is a $\,2\times 2\,$ matrix ring (over
some other ring). However, this is not the case. For instance, let
$\,R\,$ be Hamilton's quaternion division algebra generated over
$\,{\mathbb R}\,$ by $\,i,j$, with the relations $\,i^2=j^2=-1\,$
and $\,ij=-ji$. Let $\,u=a+b\,i\;(a,b\in {\mathbb R})\,$ be a
primitive $\,n$-th of unity in $\,{\mathbb R}\,[\,i\,]\cong
{\mathbb C}$. If $\,n\geq 3$, then $\,b\neq 0$, and $\,[u,j]
=[a+b\,i,j]=b\,[i,j]=2\,b\,ij\in {\rm U}(R)$. This shows that
$\,u\in R\,$ is completable. However, being a division algebra,
$\,R\,$ is not a $\,2\times 2\,$ matrix ring over any ring.
In a similar vein, if a ring $\,R\,$ has a completable element
$\,w\,$ such that $\,w^n\in \{0,w\}\,$ for some $\,n\geq 3$,
$R\,$ need not be a $\,2\times 2\,$ matrix ring either. For
instance, let $\,R={\mathbb M}_3(S)\,$ over a ring $\,S$.
For the two matrices
$$
w={\footnotesize{\begin{pmatrix}1&1&0\\0&0&1\\0&0&-1
\end{pmatrix}} } \;\,\mbox{and} \;\;
w_1={\footnotesize{\begin{pmatrix}1&0&0\\0&0&0\\1&1&0\end{pmatrix}} },
\;\mbox{we have}\;\,
[w,w_1]={\footnotesize{\begin{pmatrix}0&-1&0\\1&1&0\\-2&-2&-1
\end{pmatrix}} } \in {\rm U}(R).
$$
Thus, $\,w\,$ is completable in $\,R$, and we can check easily
that $\,w^3=w$. However, if we choose $\,S\,$ to be a division
ring (or $S={\mathbb Z}$), $R\,$ is {\it not\/} a $\,2\times 2\,$
matrix ring over any ring. The case where $\,w^n=0\,$ for some
$\,n\geq 3\,$ can be handled similarly by taking instead $\,w
=\footnotesize{\begin{pmatrix}0&1&0\\0&0&1\\0&0&0\end{pmatrix}}$
and $\,w_2=\footnotesize{\begin{pmatrix}0&0&0\\1&0&1\\1&0&0
\end{pmatrix}}$, with $\,w^3=0\,$ and $\,[w,w_2]
=\footnotesize{\begin{pmatrix}1&0&1\\1&-1&0\\0&-1&0\end{pmatrix}}
\in {\rm U}(R)$.
   
\sk

\bs\nt
{\bf Acknowledgments.} We are grateful to Professor S.\,K.~Berberian
who kindly communicated to us his detailed personal notes on
Kaplansky's ``Exercise 6'' from [Ka:~p.\,25]. These notes helped
us come to a better understanding of Kaplansky's proposed solution
of his exercise, which eventually led to the present work.  We
thank Professor G.\,M.~Bergman for his valuable contribution
[Be$_2$] toward the validation of Example 3.6, and we also thank
the referee of this paper for his/her many thoughtful comments.

\mk\nt

\sk\nt
Footnotes:

\mk\nt
{\it {\rm 2010\/} AMS Subject Classification\/}: 12E15, 15B33, 
15B36, 16E50, 16N40, 16U60.\\
{\it Keywords\/}: Units, idempotents, commutators, anti-commutators,
Bott-Duffin inverse, matrix rings, simple artinian rings, exchange
rings, commutative rings.

\bs

\bs
\nt Department of Mathematics \\
\nt Panjab University \\
\nt Chandigarh 160\,014, India

\sk\nt 
{\tt dkhurana@pu.ac.in}

\bs
\nt Department of Mathematics \\
\nt University of California \\
\nt Berkeley, CA 94720, USA   

\sk\nt 
{\tt lam@math.berkeley.edu}


\begin{thebibliography}{ABC}
  
\bibitem[Be$_1$]{Be1} G.\,M.~Bergman: {\it Modules over coproducts
of rings.}  Trans. Amer. Math. Soc. {\bf 200} (1974), 1--32.

\bibitem[Be$_2$]{Be2} G.\,M.~Bergman: Email communication,
July 3, 2017.

\bibitem[BD]{BD} R.~Bott and R.\,J.~Duffin: {\it On the algebra
of networks.} Trans. Amer. Math. Soc. {\bf 74} (1953), 99--109.

\bibitem[Co$_1$]{Cn1} P.\,M.~Cohn: {\it Some remarks on the
invariant basis property.} Topology {\bf 5} (1966), 215--228.

\bibitem[Co$_2$]{Co2} P.\,M.~Cohn: {\it Some remarks on 
projective-free rings.} Algebra Universalis {\bf 49} (2003), 159--164.

\bibitem[Co$_3$]{Co3} P.\,M.~Cohn: {\it Another criterion for a ring 
to be projective-free.} Bull.~London Math.~Soc. {\bf 37} (2005),
857--859.

\bibitem[FMP]{FMP} P.\,R.~Fuchs, C.\,J.~Maxson and G.\,F.~Pilz:
{\it On rings for which homogeneous maps are linear.}
Proc. Amer. Math. Soc. {\bf 112} (1991), 1--7.

\bibitem[Go]{Go} K.\,R.~Goodearl: {\it Von Neumann Regular 
Rings.} Second Edition, Robert E.~Krieger Publishing Company, 
Malabar, Florida, 1991.

\bibitem[GH]{GH} K.\,R.~Goodearl and D.~Handelman: {\it Simple 
self-injective rings.}  Comm.~Alg. {\bf 3} (1975), 797--834.

\bibitem[GW]{GW} H.~Grover, Z.~Wang, D.~Khurana, J.~Chen and
T.\,Y.~Lam:  {\it Sums of units in rings.} J. Alg. Appl. 
{\bf 13} (2014), 1350072, 10 pp.

\bibitem[He]{He} M.~Henriksen: {\it Two classes of rings 
generated by their units.} J.~Alg. {\bf 31} (1974), 182--193.

\bibitem[Ka]{Ka} I.~Kaplansky: {\it Rings of Operators.}
W.\,A.~Benjamin, New York, 1968.

\bibitem[Kt]{Kt} T.~Kato: {\it Perturbation Theory for Linear 
Operators.} Second Ed., Classics in Mathematics, Springer-Verlag, 
Berlin-Heidelberg-New York, 1995.

\bibitem[KL]{KL} D.~Khurana and T.\,Y.~Lam: {\it Invertible 
commutators in matrix rings.} J.~Algebra Appl. {\bf 10} 
(2011), 51--71.  

\bibitem[KLS]{KLS} D.~Khurana, T.\,Y.~Lam and N.~Shomron: 
{\it A quantum-trace determinantal formula for matrix commutators,
and applications.} Lin. Alg. Appl. {\bf 436} (2012), 2380--2397.

\bibitem[KR$_1$]{KR1}  J.\,J.~Koliha and V.~Rako\v{c}evi\'c: 
{\it Invertibility of the sum of idempotents.} 
Lin. Multilin.~Algebra {\bf 50} (2002), 285--292.

\bibitem[KR$_2$]{KR2}  J.\,J.~Koliha and V.~Rako\v{c}evi\'c: 
{\it Invertibility of the difference of idempotents.} 
Lin. Multilin.~Algebra {\bf 51} (2003), 97--110.

\bibitem[La$_1$]{La1} T.\,Y.~Lam: {\it A First Course in 
Noncommutative Rings.}  Second Edition, Graduate Texts in Math., 
Vol.~{\bf 131}, Springer-Verlag, Berlin-Heidelberg-New York, 2001.

\bibitem[La$_2$]{La2}  T.\,Y.~Lam: {\it Exercises in Classical 
Ring Theory.}  Second Edition, Problem Books in Mathematics, 
Springer-Verlag, Berlin-Heidelberg-New York, 2003.

\bibitem[La$_3$]{La3} T.~Y.~Lam: {\it Lectures on Modules and 
Rings.} Graduate Texts in Math., Vol.~{\bf 189}, Springer-Verlag, 
Berlin-Heidelberg-New York, 1999. 

\bibitem[La$_4$]{La4}  T.\,Y.~Lam: {\it Exercises in Modules
and Rings.} Problem Books in Mathematics, Springer-Verlag, 
Berlin-Heidelberg-New York, 2007.

\bibitem[Le]{Lea} W.\,G.~Leavitt: {\it Rings without invariant 
basis number.} Proc. Amer. Math. Soc. {\bf 8} (1957), 322--328.

\bibitem[Ni]{Ni} W.\,K.~Nicholson: {\it Lifting idempotents and 
exchange rings.} Trans. Amer. Math. Soc. {\bf 229} (1977), 269--278.
  
\bibitem[Sa]{Sa} E.\,M.~de S\'a: {\it The rank of a difference
of similar matrices.} Portugal. Math. {\bf 46} (1989), 177--187.
  
\bibitem[SG]{SG} G.~Song and X.~Guo: {\it Diagonability
of idempotent matrices over noncommutative rings.}
Lin. Multilin. Algebra {\bf 297} (1999), 1--7.

\bibitem[Wa]{Wa} R.\,B.~Warfield: {\it Exchange rings and 
decompositions of modules.} Math.~Ann. {\bf 199} (1972), 31--36.

\end{thebibliography}
\end{document}